\documentclass{amsart}
\usepackage{amssymb}
\DeclareFontFamily{OMS}{cmsy}{%
\fontdimen16\font=3pt
\fontdimen17\font=3pt}
\makeatletter
\renewcommand{\subsection}{\@startsection{subsection}{2}{\z@}%
{\baselineskip}{0.5\baselineskip}{\bfseries}}
\makeatother
\def\dj{d\kern-.30em\raise1.25ex\vbox{\hrule width .3em height .03em}}
\def\Dj{D\rlap{\kern-.70em\raise0.75ex
\vbox{\hrule width .3em height .03em}}}
\def\k{\kappa}
\def\cal{\mathcal}
\def\Bbb{\mathbb}
\def\frak{\mathfrak}
\newtheorem{pro}{Proposition}
\newtheorem{lem}[pro]{Lemma}
\def\hor{\frak{hor}}
\def\ver{\frak{ver}}
\def\con{\mathrm{con}}
\def\adj{\varpi}
\def\sstar{{\raise0.2ex\hbox{$\scriptstyle\star$}}}
\def\HZ{H\!Z}
\def\3{\vartriangle}

\def\S{\Sigma}
\def\J{J_{\!\sigma}}
\def\IS{I(\Sigma)}
\def\IG{I^\sstar}
\def\bla#1{$(${\it #1\/{}}$)$}
\def\k{\kappa}
\def\im{\mathrm{im}}
\def\WG{\Omega}
\def\WGi{\daleth}
\def\fWG{\widetilde{\adj}}
\def\hh{\frak{h}}
\def\rWG{\adj_\wedge}
\def\WGJ{\cal{K}}
\def\WGr{\Omega_*}
\def\WGri{\beth}
\newcommand{\grten}{\mathbin{\widehat{\otimes}}}
\newcommand{\vh}{\frak{vh}}
\newcommand{\ad}{\mathrm{ad}}
\newcommand{\id}{\mathrm{id}}
\newcommand{\inv}{{i\!\hspace{0.8pt}n\!\hspace{0.6pt}v}}
\def\e{\epsilon}
\newcommand{\Sum}{{\displaystyle\sum}}
\def\Gr{\frak{gr}}
\newcommand{\kI}{\Rsh}
\newcommand{\kW}{\amalg}
\begin{document}
\title[Quantum Principal Bundles]{Quantum Principal Bundles $\&$ Their\\
Characteristic Classes}
\author{Mi\'co {\Dj}ur{\Dj}evi\'c}
\address{Instituto de Matematicas, UNAM, Area de la Investigacion
Cientifica, Circuito Exterior, Ciudad Universitaria, M\'exico DF, cp
04510, MEXICO}
\thanks{Lectures presented at Workshop on Quantum and Classical Gauge
Theories, Stefan Banach International Mathematical Center, Warsaw,
Poland, May 1995.}
\begin{abstract}
A general theory of characteristic classes of quantum principal bundles
is sketched, incorporating basic ideas of classical Weil theory into the
conceptual framework of non-commutative differential geometry.
A purely cohomological interpretation of the Weil
homomorphism is given, together with a geometrical interpretation via
quantum invariant polynomials. A natural spectral
sequence is described. Some quantum phenomena appearing in
the formalism are discussed.
\end{abstract}
\maketitle
\section{Introduction}

In this letter we are going to present the basic structural elements of a
noncommutative geometric \cite{C1,C2} generalization of classical
Weil theory \cite{KN} of characteristic classes. All considerations are
logically based on a general theory of quantum principal
bundles, developed in \cite{D1,D2}. This paper does not contain the
proofs. A detailed exposition of the theory of quantum characteristic
classes can be found in~\cite{D3}.

A general construction of characteristic classes
is presented in Section~3. If the calculus on the bundle admits certain
special connections (satisfying properties called regularity and
multiplicativity), the classical Weil construction of characteristic classes
(via invariant polinomials) can be directly incorporated \cite{D2}
into the quantum context. We shall first consider this particular case.
Then we present a different construction of
characteristic classes, in the framework of which a natural cohomological
interpretation of the domain of the Weil homomorphism is given. This
construction works
for arbitrary quantum principal bundles and connections.
Section 4 is devoted to the study of a natural spectral sequence,
which can be associated to an arbitrary quantum principal bundle.

In the next section basic elements of the differential calculus over
quantum principal bundles and the theory of connections are collected.
The last section contains concluding examples. A particular
attention is given to the purely quantum phenomena appearing in the
formalism. Finally, the Appendix contains the list of quantum group
symbols figuring in the main text.

\smallskip
\section{Geometrical Background}
\renewcommand{\thepage}{}

Let us consider a compact matrix \cite{W1} quantum group $G$, represented
by a Hopf *-algebra $\cal{A}$, interpreted as consisting
of polynomial functions on $G$. Let
$\phi,\e$ and $\k$ be the coproduct, counit and the antipode map.

Let $M$ be a quantum space, represented by a *-algebra
$\cal{V}$. Let $P=(\cal{B},i,F)$ be a quantum principal $G$-bundle over
$M$. Here $\cal{B}$  is a *-algebra representing $P$ as a quantum space,
while $i\colon\cal{V}\rightarrow\cal{B}$ and
$F\colon\cal{B}\rightarrow\cal{B}\otimes\cal{A}$ are unital
*-homomorphisms playing the role of the dualized projection of $P$
on $M$ and the right action of $G$ respectively.

Let us assume that the complete differential calculus on
$G$ is described by
\cite{D1} the universal envelope $\Gamma^\wedge$ of a
first-order \cite{W2} bicovariant *-calculus $\Gamma$.
Let $\Omega(P)$ be a graded-differential *-algebra
representing the calculus on $P$. By definition \cite{D2},
this means that $\Omega^0(P)=\cal{B}$, and
that $\cal{B}$ generates the differential algebra $\Omega(P)$.
Furthermore, we require that the right action $F$ is extendible to
a graded-differential homomorphism $\widehat{F}\colon\Omega(P)
\rightarrow\Omega(P)\grten\Gamma^\wedge$. This extension is unique and
hermitian.

\renewcommand{\thepage}{\arabic{page}}
The formula
$F^\wedge=(\id\otimes p_*)\widehat{F}$,
where $p_*\colon\Gamma^\wedge\rightarrow\cal{A}$ is the projection map,
defines the right action
$F^\wedge\colon\Omega(P)\rightarrow\Omega(P)\otimes\cal{A}$
of $G$ on $\Omega(P)$.

Let $\omega\colon\Gamma_{\inv}\rightarrow\Omega(P)$
be a connection on $P$. This means that $\omega$ is a hermitian intertwiner
between the adjoint action $\adj\colon\Gamma_{\inv}\rightarrow\Gamma_{\inv}
\otimes\cal{A}$ and $F^\wedge$, such that
$$\pi_v\omega(\vartheta)=1\otimes\vartheta$$
where $\pi_v\colon\Omega(P)\rightarrow\ver(P)$ is the 
verticalization homomorphism, and $\ver(P)$ is the 
graded-differential *-algebra representing `verticalized' differential 
forms on the bundle (we can write $\ver(P)=\cal{B}
\otimes\Gamma_{\inv}^\wedge$, at the level of graded vector
spaces). The set $\con(P)$ of all connections on $P$
is non-empty. It is naturally a real affine space.
Equivalently, we can say that connections on $P$ are hermitian
first-order linear maps $\omega\colon\Gamma_{\inv}\rightarrow\Omega(P)$
satisfying the identity
$$
\widehat{F}\omega(\vartheta)=\sum_k\omega(\vartheta_k)\otimes
c_k+1\otimes\vartheta,
$$
where $\adj(\vartheta)=\Sum_k\vartheta_k\otimes c_k$.

The elements of a graded *-subalgebra
$$\hor(P)=\widehat{F}^{-1}\Bigl\{\Omega(P)\otimes\cal{A}\Bigr\},$$
are naturally interpretable as horizontal forms on $P$. The above algebra
is $F^\wedge$-invariant. In the framework of this interpretation, it is
natural to define differential forms on $M$ as $F^\wedge$-invariant
elements of $\hor(P)$. These elements form a graded-differential
*-subalgebra $\Omega(M)$ of $\Omega(P)$. Equivalently, $\Omega(M)$ is
the $\widehat{F}$-fixed-point subalgebra of $\Omega(P)$.

We say that $\omega$ is {\it regular} iff
$$
\omega(\vartheta)\varphi=(-1)^{\partial\varphi}\sum_k\varphi_k\omega(\vartheta
\circ c_k)
$$
for each $\varphi\in\hor(P)$, where $\Sum_k\varphi_k\otimes
c_k=F^\wedge(\varphi)$. Let us observe that regular connections
graded-commute with the elements of $\Omega(M)$.

We say that $\omega$ is {\it multiplicative} iff it is extendible
to a unital homomorphism $\omega^\wedge\colon
\Gamma^\wedge_{\inv}\rightarrow\Omega(P)$. This extension is unique, and
hermitian.

It turns out that regular connections are all multiplicative or
not at the same time. If the calculus on the bundle is such that
regular connections are not multiplicative then
it is possible to `renormalize' it, by factorizing $\Omega(P)$
through an appropriate graded-differential *-ideal, which measures a
lack of multiplicativity of regular connections. Such a factorization
does not change the first-order calculus. In terms of the
renormalized calculus, regular connections are multiplicative. In what
follows it will be assumed that regularity implies multiplicativity.

The above definition of the map $\omega^\wedge$ can be extended to
non-multiplicative connections, with the help of the formula
$$\omega^\wedge=\omega^\otimes\iota.$$
Here $\iota\colon\Gamma_{\inv}^\wedge\rightarrow\Gamma_{\inv}^\otimes$
is a fixed grade-preserving hermitian section of the
factorization map, which intertwines the adjoint actions
of $G$, while $\omega^\otimes\colon\Gamma_{\inv}^\otimes\rightarrow\Omega(P)$
is the unital multiplicative extension of $\omega$.

The formulas
\begin{align*}
(\psi\otimes\eta)(\varphi\otimes\vartheta)&=
(-1)^{\partial\eta\partial\varphi}\sum_k\psi\varphi_k\otimes
(\eta\circ c_k)\vartheta\\
(\varphi\otimes\vartheta)^*&=\sum_k\varphi_k^*\otimes(\vartheta^*\circ
c_k^*)
\end{align*}
introduce a graded *-algebra structure in the graded vector space
$$\vh(P)=\hor(P)\otimes\Gamma_{\inv}^\wedge.$$

Let $m_\omega\colon\vh(P)\rightarrow\Omega(P)$ be a linear map defined by
$$ m_\omega(\varphi\otimes\vartheta)=\varphi\omega^\wedge(\vartheta).$$
It turns out that $m_\omega$ is bijective.
It intertwines the corresponding actions of $G$.
Moreover, if $\omega$ is regular and multiplicative then $m_\omega$ is a
*-algebra isomorphism.

The quantum counterpart of the horizontal projection
operator can be defined by
$$h_\omega=(\id\otimes p_*)m_\omega^{-1}.$$
It projects $\Omega(P)$ onto $\hor(P)$.
Following classical differential geometry, the covariant
derivative operator is naturally defined by
$$D_\omega=h_\omega d.$$
Finally, the curvature operator
$R_\omega\colon\Gamma_{\inv}\rightarrow \Omega(P)$
can be defined as the composition
$$R_\omega=D_\omega\omega.$$
The identity
$$R_\omega=d\omega-\langle\omega,\omega\rangle$$
is the counterpart of the classical structure equation.
Here $\langle\,\rangle$ are the
brackets naturally associated to the embedded differential
$\delta\colon\Gamma_{\inv}\rightarrow\Gamma_{\inv}\otimes \Gamma_{\inv}$,
given by composing $d$ and the embedding $\iota$.
Let us observe that the curvature operator
$R_\omega$ implicitly depends on the map $\delta$.

Let $\IG\subseteq\Gamma_{\inv}^\otimes$ be the graded *-subalgebra
consisting of elements invariant under the adjoint action $\adj\colon
\Gamma^\otimes_{\inv}\rightarrow\Gamma_{\inv}^\otimes\otimes\cal{A}$. The
formula
$$
W^\omega(\vartheta)=R_\omega^\otimes(\vartheta)
$$
defines a *-homomorphism $W^\omega\colon\IG\rightarrow\Omega(M)$.
Here $R_\omega^\otimes\colon\Gamma_{\inv}^\otimes\rightarrow\Omega(P)$ is
the corresponding unital multiplicative extension.

We shall denote by $Z(M)$ the graded centre of
$\Omega(M)$. It is a graded-differential *-subalgebra of $\Omega(M)$.

If $\omega$ is regular then the following identity holds
$$ R_\omega(\vartheta)\varphi=\sum_k\varphi_kR_\omega(\vartheta\circ
c_k),$$
for each $\varphi\in\hor(P)$.
In particular, the curvature $R_\omega$ commutes with all elements of
$\Omega(M)$, and hence the image of $W^\omega$ is contained in $Z(M)$.
The covariant derivative of regular connections is hermitian, and
satisfies the graded Leibniz rule.

\section{Quantum Characteristic Classes}

\subsection{The Regular Case}

Let us assume that the bundle admits regular (and multiplicative)
connections. A large class of examples of such bundles is given by
non-commutative frame structures \cite{D4}. In this special case there
exists a canonical regular connection $\omega$, corresponding to
the Levi-Civita connection in classical geometry.

Let us assume that $\omega$ is a regular connection on a quantum
principal bundle $P$ with the calculus $\Omega(P)$.
Applying essentially the same transformations as in the
classical case \cite{KN}, it follows that
$W^{\omega}(\vartheta)$ is closed, for each $\vartheta \in \IG$.

\begin{lem}
The cohomological
class of $W^{\omega}(\vartheta)$ in $Z(M)$ is
independent of the choice of a regular connection $\omega$.
\end{lem}

In other words, there exists the intrinsic unital *-homomorphism
$$W\colon\IG\rightarrow \HZ(M)\qquad
W(\vartheta)=[W^{\omega}(\vartheta)],$$
where $\HZ(M)$ is the corresponding cohomology algebra. This map is a
quantum counterpart of the Weil homomorphism.

The Weil homomorphism is further factorizable through
the ideal $\J$  generated  by the space
$\im(I-\sigma)\subseteq \Gamma_{\inv}^{\otimes 2}$. This follows
from the commutation identity
$$ R_\omega(\eta)R_\omega(\vartheta)=\sum_k
R_\omega(\vartheta_k)R_\omega(\eta\circ c_k) $$
where $\adj(\vartheta)=\Sum_k\vartheta_k\otimes c_k$.
The elements of the factoralgebra
$\S =\Gamma_{\inv}^\otimes/\J$
are quantum counterparts of the polynoms over the Lie  algebra of the
structure group.
The  adjoint action $\adj$ is naturally projectable to
$\adj_{\S}\colon\S\rightarrow\S\otimes\cal{A}$.
Let us denote by $\IS\subseteq \S$
the  subalgebra  of  elements  invariant  under $\adj_{\S}$
(invariant polynomials). We have
$\IS=\IG\big/\bigl(\IG\cap \J\bigr)$, and
from the commutation relations
defining $\S$ it follows that $\IS$ is a {\it central}
subalgebra of $\S$.

\subsection{The General Case}

In generalising the theory of characteristic classes
to the level of arbitrary bundles, we shall follow {\it the idea of
universality}. Algebraic expressions generating characteristic
classes should be the same for all bundles. The following observation
is the starting point.

\begin{pro}
Let us assume that $\omega$ is a regular and multiplicative connection
on $P$. Then the image of the map $W^\omega$ consists of differential
forms on $M$ which are expressible in terms of $\omega$ and $d\omega$.
\end{pro}

Let $\WG$ be a differential algebra generated
by the first-order space $\Gamma_\inv$ with the
only relation $d(1)=0$. The
*-involution on $\Gamma_\inv$ naturally extends to $\WG$, so that
$d\colon\WG\rightarrow\WG$ is a hermitian map. The cohomology of $\WG$
is trivial--we have $ H(\WG)=\Bbb{C}$.

There exists the unique graded-differential homomorphism
$\fWG\colon\WG\rightarrow\WG\grten\Gamma^\wedge$ such that
$$
\fWG(\vartheta)=\adj(\vartheta)+1\otimes\vartheta
$$
for each $\vartheta\in\Gamma_{\inv}$. Furthermore this map is hermitian and
$$
(\fWG\otimes\id)\fWG=(\id\otimes\widehat{\phi})\fWG.
$$

It is also possible to introduce a natural right action $\rWG\colon
\WG\rightarrow\WG\otimes\cal{A}$, extending the adjoint action map
$\adj$. Explicitly,
$$\rWG=(\id\otimes p_*)\fWG. $$
Let $\WGi\subseteq\WG$ be the $\fWG$-fixed-point graded differential
*-subalgebra.

Let $P=(\cal{B},i,F)$ be a quantum principal $G$-bundle over $M$,
endowed with a calculus $\Omega(P)$. Let $\omega$ be an arbitrary
connection on $P$.
\begin{pro}
\bla{i} There exists the unique homomorphism
$\widehat{\omega}\colon\WG\rightarrow\Omega(P)$ of differential
algebras extending $\omega$.
This map is hermitian. Moreover, it intertwines $\fWG$ and
$\widehat{F}$. In particular, it follows that
$\widehat{\omega}(\WGi)\subseteq\Omega(M)$.

\smallskip
\bla{ii} The induced cohomology map $W\colon H(\WGi)\rightarrow H(M)$ is
independent of the choice of $\omega$.
\end{pro}

The constructed map $W$ is a counerpart of the Weil homomorphism, at the
level of general quantum principal bundles. Characteristic classes are
therefore labeled by the elements of $H(\WGi)$.

The algebra $\WG$ possesses various properties characteristic to
differential algebras describing the calculus on quantum principal
bundles. In particular, it is possible to introduce a natural decomposition
$$
\WG\leftrightarrow\hh(\WG)\otimes\Gamma_{\inv}^\wedge=\vh(\WG)\qquad
\varphi\iota(\vartheta)\leftrightarrow\varphi\otimes\vartheta
$$
where $\hh(\WG)\subseteq\WG$ is a graded $*$-subalgebra describing
`horizontal elements', defined by
$$\hh(\WG)=\fWG^{-1}(\WG\otimes\cal{A}).$$
It follows that $\rWG[\hh(\WG)]\subseteq\hh(\WG)\otimes\cal{A}$, in other
words $\hh(\WG)$ is $\rWG$-invariant. The algebra $\WGi$ is the
$\rWG$-fixed point subalgebra of $\hh(\WG)$.

The first step in computing characteristic classes is to find the cocycles
of $\WGi$. The following is a prescription of constructing these cocycles.
Every cocycle $w\in\WGi^+$ is of the form $w=d\varphi$, where
$\varphi$ is some $\rWG$-invariant element of $\WG$. Then we have the
equivalence
$$ w\in\WGi\iff d\bigl\{\fWG(\varphi)-\varphi\otimes 1\bigr\}=0.$$

\begin{lem} Let $\cal{C}\subseteq\WG\otimes\Gamma^{\wedge+}$ be the a
subcomplex spanned by elements of the form
$c=\fWG(\varphi)-\varphi\otimes 1$, where $\varphi\in\WG$ is
$\rWG$-invariant. Then the following natural correspondence holds
$$
H^n(\WGi)\leftrightarrow H^{n-1}(\cal{C}).
$$
\end{lem}

In particular, for every quantum principal bundle $P$
the cocycles representing characteristic classes are exact,
as classes on the bundle, with
$\widehat{\omega}(w)=d\widehat{\omega}(\varphi)$. Therefore, the above
introduced elements $\varphi$ play the role of {\it universal} Chern-Simons
forms.

Another possible approach in computing quantum characteristic classes
consists in finding explicitly horizontal forms, and computing invariant
elements ($\Leftrightarrow$ the algebra $\WGi$).
The algebra $\hh(\WG)$ is invariant under the action of the operator
$$
\ell_\star(\vartheta,\varphi)=\vartheta\varphi-(-1)^{\partial
\varphi}\sum_k\varphi_k(\vartheta\circ c_k),
$$
where $\vartheta\in\Gamma_{\inv}$ and $\Sum_k\varphi_k\otimes c_k=
\rWG(\varphi)$.

\begin{pro} The algebra $\hh(\WG)$
is the minimal $\ell_\star$-invariant subalgebra of $\WG$ containing
the elements $R(\vartheta)=d\vartheta-\delta(\vartheta)$ and the elements
from $S_{\inv}^{\wedge 2}$.
\end{pro}

Let us now return to the structures admiting regular and multiplicative
connections. If we work only with regular connections, then the
cohomological construction of the Weil homomorphism should be appropriately
refined, by factorizing
$\WG$ through the ideal which takes into account the regularity
property.

Let $\WGJ\subseteq\WG$ be the ideal generated by elements from
$S_{\inv}^{\wedge2}$, and the elements of the form
\begin{gather*}
j_3(\eta,\vartheta)=\ell_\star(\eta,R(\vartheta))=\eta
R(\vartheta)-\sum_k R(\vartheta_k)(\eta\circ c_k)\\
j_4(\eta,\vartheta)=R(\eta)R(\vartheta)-\sum_k R(\vartheta_k)R
(\eta\circ c_k),
\end{gather*}
where $\eta,\vartheta\in\Gamma_{\inv}$ and
$\adj(\vartheta)=\Sum_k\vartheta_k\otimes c_k$.

It follows that $\Gamma_{\inv}^\wedge$ is a subalgebra of
$\WGr=\WG/\WGJ$, in a natural manner. We have
$$
\fWG R(\vartheta)=\sum_kR(\vartheta_k)\otimes c_k\qquad
R(\vartheta)^*=R(\vartheta^*).
$$
The ideal
$\WGJ$ is invariant under the actions of $\fWG,*$ and $d$.

The maps $\fWG,d,*,\rWG$ are hence projectable to $\WGr$. Obviously,
projected maps (that will be denoted by the same symbols)
are in the same algebraic relations as the original ones.
Let us introduce the horizontal part
$\hh(\WGr)=\fWG^{-1}(\WGr\otimes\cal{A})$ of $\WGr$,
which is a $\rWG$-invariant *-subalgebra of $\WGr$. Applying a similar
reasoning as in \cite{D2} it follows that
\begin{lem} \bla{i} The product map in $\WGr$ induces a graded vector
space isomorphism
$$
\hh(\WGr)\otimes\Gamma_{\inv}^\wedge\leftrightarrow\WGr.
$$

\smallskip
\bla{ii} The map $R\colon\Gamma_{\inv}\rightarrow\WGr$ can be uniquely extended
to a unital *-homomorphism $R\colon\S\rightarrow\WGr$. The extended
$R$ maps isomorphically $\S$ onto $\hh(\WGr)$. Moreover,
$R$ intertwines $\adj_{\S}$ and $\rWG$.
\end{lem}

Let us formulate a regular counterpart of Proposition~3.
Let us consider a quantum principal $G$-bundle $P=(\cal{B},i,F)$
over $M$, with a calculus $\Omega(P)$ admiting regular and
multiplicative connections. Let $\omega$ be an arbitrary regular
connection.

\begin{pro}
\bla{i} There exists the unique homomorphism
$\widehat{\omega}\colon\WGr\rightarrow\Omega(P)$ of differential
algebras extending the connection $\omega$.
The map $\widehat{\omega}$ is hermitian, and
$$ (\widehat{\omega}\otimes\id)\fWG=\widehat{F}\widehat{\omega}. $$
In particular, $\widehat{\omega}(\WGri)\subseteq\Omega(M)$.

\smallskip
\bla{ii} The following identities hold
$$
\widehat{\omega}R=R_\omega\qquad
\widehat{\omega}D=D_\omega\widehat{\omega}
$$
where $D\colon\WGr\rightarrow\WGr$ is a first-order antiderivation specified
by $$DR(\vartheta)=0\qquad D\vartheta=R(\vartheta).$$

\smallskip
\bla{iii} The induced cohomology map $W\colon H(\WGri)\rightarrow H(M)$ is
independent of the choice of $\omega$.

\smallskip
\bla{iv} We have $H(\WGri)=\WGri=\IS$, in a natural manner. This gives
a connection with the definition of the Weil homomorphism given in the
previous subsection.
\end{pro}

The above introduced map $D$ is called the {\it
universal covariant derivative}. We have $D^2=0$.

Let $d_{v\!h}\colon\WGr\rightarrow\WGr$ be `the universal'
vertical differential. By definition, this map is acting in the
following way
$$
d_{v\!h}(\varphi\otimes\vartheta)=\sum_k\varphi_k\otimes
\pi(c_k)\vartheta+\varphi\otimes d(\vartheta),
$$
where $\rWG(\varphi)=\Sum_k\varphi_k\otimes c_k$.
\begin{lem}
The following identities hold
\begin{align*}
Dd_{v\!h}+d_{v\!h}D&=0\\
D+d_{v\!h}&=d.
\end{align*}
\end{lem}

Let us observe that the map $D$ acts skew-diagonally, with respect to a
natural bigrading in $\WGr$. This implies that the cohomology algebra
$H_{\!D}(\WGr)$ is naturally bigraded, too.

\section{The Spectral Sequence}

Let us consider a quantum principal bundle $P$ endowed with a
differential structure $\Omega(P)$.

For each $k\geq 0$ let $\Omega_k(P)\subseteq\Omega(P)$ be the space
consisting of elements having the `vertical order' less or equal $k$.
In other words, $$\Omega_k(P)=\widehat{F}^{-1}\bigl(\Omega(P)\otimes
\Gamma^\wedge_k\bigr)$$
where $\Gamma_k^\wedge$ consists of forms having degrees not exceeding
$k$. These spaces form a filtration of $\Omega(P)$, and the
following compatibility properties hold:
$$\Omega_k(P)^*=\Omega_k(P)\qquad d\Omega_k(P)\subseteq\Omega_{k+1}(P)
\qquad\Omega_k(P)=\sideset{}{^\oplus}\sum_{j\geq 0}\Omega_k^j(P).$$
Let us consider a graded-differential *-algebra
$$\kW(P)=\sideset{}{^\oplus}\sum_{k\geq 0}\Omega_k(P), $$
where the grading is given by numbers $k$, and the differential
*-structure is induced from $\Omega(P)$. Let $\kI\colon\kW(P)\rightarrow
\kW(P)$ be the first-order map induced by the inclusions
$\Omega_k(P)\subseteq\Omega_{k+1}(P)$. By definition, this is a
monomorphism of differential *-algebras. Let $\Gr(P)$ be the
graded-differential *-algebra associated to the introduced filtration.
In other words, we have a short exact sequence
$$
0\xrightarrow{}\kW(P) \xrightarrow{\kI}\kW(P)\xrightarrow{}\Gr(P)
\xrightarrow{}0
$$
of differential *-algebras.

The space $\Omega_k(P)$ is linearly spanned by elements of the form
$$w=\varphi\omega(\vartheta_1)\dots\omega(\vartheta_j)$$ where
$\varphi\in\hor(P)$ and $j\leq k$, while $\omega$ is an arbitrary
connection. The algebra $\Gr(P)$ is invariantly isomorphic to the
algebra $\vh(P)$ of `vertically-horizontally' decomposed forms.
Explicitly, the isomorphism $\vh(P)\leftrightarrow\Gr(P)$ is
given by
$$ \bigl(\varphi\otimes(\vartheta_1\dots\vartheta_k)\bigr)
\leftrightarrow\bigl[\varphi\omega(\vartheta_1)\dots\omega(\vartheta_k)
+\Omega_{k-1}(P)\bigr].$$
We shall assume that the two algebras are identified,
with the help of the above isomorphism. We shall also assume that the
calculus $\Gamma$ is such that only scalar elements of
$\cal{A}$ are anihilated by the differential map.

The factor-differential $d_{v\!h}\colon\vh(P)\rightarrow\vh(P)$
is given by $$d_{v\!h}(\varphi\otimes\vartheta)=
(-1)^{\partial\varphi}\sum_k\varphi_k\otimes\pi(c_k)\vartheta+
(-1)^{\partial\varphi}\varphi\otimes d(\vartheta)$$
where $\Sum_k\varphi_k\otimes c_k=F^\wedge(\varphi)$.

It turns out that
$$
H\bigl[\vh(P)\bigr]=\Omega(M)\grten H(\Gamma^\wedge_{\inv}),
$$
in a natural manner.

Let $E(P)=\Bigl\{E_r(P)\vert r\in\Bbb{N}\Bigr\}$ be the spectral
sequence associated to the introduced short exact sequence.
The introduced filtration of $\Omega(P)$ induces a filtration of the
*-algebra $H(P)$ of cohomology classes. We have
$$H_k(P)= \Sum_j^\oplus H^j_k(P).$$
Applying general theory \cite{bot}, it
follows that the introduced spectral sequence is convergent, and that
$E_\infty(P)$ coincides with the graded *-algebra associated to the
filtered $H(P)$.

By construction $E_1(P)=H\bigl[\vh(P)\bigr]$. Furthermore, it turns out that
the differential $d_1$ is given by
$$
d_1(w\otimes[\vartheta])=dw\otimes[\vartheta].
$$
In particular,
$$
E_2(P)=H(M)\grten H(\Gamma_{\inv}^\wedge).
$$

A similar considerations can be applied to algebras $\WG$ and $\WGr$, in
particular we can associate natural spectral sequences to these
algebras. The spectral sequence $E(\WG)$ converges to the trivial cohomology
$H(\WG)=\Bbb{C}$. Therefore the spaces $H^k(\WGi)$ are naturally
filtered. 

The trivial convergence information is insufficient to compute the
cohomology algebra $H(\WGi)$. However, in various interesting special cases
the spectral sequence degenerates (as in classical geometry), and the
triviality property is sufficient to determine all cohomology classes.

\section{Concluding Remarks}

We have assumed that the higher-order calculus on $G$ is
described by the universal envelope $\Gamma^\wedge$. This
corresponds to the {\it maximal} solution. All the constructions with
differential forms on quantum principal bundles can be performed
dealing with the bicovariant \cite{W2} exterior algebra
$\Gamma^\vee$, instead of $\Gamma^\wedge$. This corresponds to
the {\it minimal} appropriate higher-order calculus, the essential property is
that the coproduct map is extendible to a homomorphism
$\phi^\vee\colon\Gamma^\vee\rightarrow\Gamma^\vee{\grten}
\Gamma^\vee$ of graded-differential *-algebras.

Moreover, all the constructions can be performed also for
`intermediate' higher-order calculi, described by
higher-order graded-differential *-ideals
$S_\star\subseteq\Gamma^\wedge$ satisfying
$$ \widehat{\phi}(S_\star)\subseteq S_\star\grten \Gamma^\wedge +
\Gamma^\wedge\grten S_\star.$$

Changing the higher-order calculus over $G$ directly influences the
algebra of higher-order horizontal forms, and in particular the
higher-order part of $\Omega(M)$.

The cohomological formulations of the
two levels of the theory (regular and general) are essentially the same.
However, it is interesting to observe that
the corresponding universal characteristic classes radically differ, if
the calculus on the group is sufficiently `irregular'.

For example, at the level of general bundles there exist generally
nontrivial classes in {\it odd} dimensions, in contrast to the regular
case where all classes are expressible in terms of the curvature map.
Another important difference between two levels is concerning the
existence of the Chern-Simons forms. At the general level, all
characteristic classes vanish, as cohomology classes of the bundle. In
contrast to this, generally there exist regular characteristic
classes, nontrivial
as classes on the bundle. In particular, such classes have no
analogs at the level of general differential structures. On the other
hand, under certain regularity assumptions on the calculus over
the structure group, universal classes for the regular and the general
case coincide. This essentially simplifies the work with bundles without
regular connections.

The higher-order calculus on $G$ can be always maximally addopted, from
the point of view of the appearance of $3$-dimensional characteristic
classes. Let us consider this point in more details.

In dimension $2$, Chern-Simons forms are labeled by closed
elements of $\Gamma_{\inv}^{\wedge 2}$, invariant under the actions of
$\adj$ and $\sigma$. Let
$S_3\subseteq\Gamma^\wedge$ be the ideal generated by the elements of
the form $d(\psi)$, where $\psi\in\Gamma_{\inv}^{\wedge 2}$
satisfy
$$\sigma(\psi)=\psi\qquad\adj(\psi)=\psi\otimes 1.$$

\begin{lem} The space $S_3$ is a graded-differential *-ideal in
$\Gamma^\wedge$. Moreover, $$\widehat{\phi}(S_3)\subseteq
S_3\grten\Gamma^\wedge+ \Gamma^\wedge\grten S_3.$$ \end{lem}

Let us pass to the factor-calculus $\Gamma^{\3}=\Gamma^\wedge/S_3$, in
the framework of which the elements $\psi$ are closed. It follows that
\begin{lem}
Relative to the factorized calculus
$$
H^3(\WGi)=\Bigl\{\text{The elements }\psi\in
\Gamma_{\inv}^{\wedge 2}\Bigr\}\Bigm/d\Bigl\{\adj\text{-invariant
elements of }\Gamma_{\inv}\Bigr\}.
$$
\end{lem}

\bigskip
\centerline{\scshape Appendix: Quantum Group Symbols}

\bigskip
This table gives a brief explanation of quantum group symbols
used in the main text.

\medskip
\def\ika{\noalign{\hrule}}
\def\vv{\vphantom{\Bigl\{}}
\def\pom{\hskip 5pt}
\vbox{\tabskip=0pt\offinterlineskip\hrule
\halign{&\vrule#&\strut\pom#\pom\hfill\cr
height2pt&\omit&&\omit&&\omit&\cr
&{\it Symbol}\hfill && {\it Description}\hfill &&
{\it Comments} \hfill&\cr
height2pt&\omit&&\omit&&\omit&\cr
\noalign{\hrule}
&$\vv\Gamma\qquad d\colon\cal{A}\rightarrow\Gamma$&&\multispan3
\!\!\pom A bicovariant first-order
*-calculus over $G$ \pom\hfill &\cr\ika
&$\Gamma^{\wedge}$&& The universal differential  &&
$\Gamma^\wedge=\Gamma^\otimes/S^\wedge$&\cr
&$d\colon\Gamma^\wedge
\rightarrow\Gamma^\wedge$ && envelope of $\Gamma$
&&$S^\wedge$--the corresponding ideal &\cr\ika
&$\Gamma^\vee$&& The braided exterior&&
$\Gamma^\vee=\Gamma^\otimes/\ker(A)$&\cr
&$d\colon\Gamma^\vee\rightarrow\Gamma^\vee$ && algebra over $\Gamma$
&&$A$--the braied antisymetrizer &\cr\ika
&$\vv\sigma\colon\Gamma_{\inv}^{\otimes 2}\rightarrow\Gamma_{\inv}^{\otimes
2}$ &&\multispan3\pom The intrinsic braid operator, factorizable to
$\Gamma_{\inv}^{\wedge 2}$\pom\hfill&\cr\ika
&$\vv\Gamma^\otimes$&&\multispan3\!\! \pom The tensor bundle algebra
over $\Gamma$\pom\hfill &\cr\ika
&$\vv\widehat{\phi}\colon\Gamma^{\wedge}\rightarrow\Gamma^{\wedge}
\grten\Gamma^{\wedge}$&&\multispan3\!\! \pom \strut The differential extension
of the coproduct map $\phi$\pom\hfill&\cr\ika
&$\vv\Gamma_{\inv}^\sstar,\pom\star\in\{\wedge,\vee,\otimes\}$&&
The left-invariant  part&&
$\Gamma_{\inv}^\wedge=\Gamma_{\inv}^\otimes/S^\wedge_{\inv}$&\cr\ika
&$\vv\pi\colon\cal{A}\rightarrow\Gamma_{\inv}$&& A natural surjection map &&
$\pi(a)=\k(a^{(1)})da^{(2)}$&\cr\ika
&$\circ\colon\Gamma_{\inv}^\sstar\otimes\cal{A}
\rightarrow\Gamma^\sstar_{\inv}$ && A natural right-module&&
$\vartheta\circ  a=\k(a^{(1)})\vartheta a^{(2)}$&\cr
& && structure map && &\cr\ika
&$\vv\ad\colon\cal{A}\rightarrow\cal{A}
\otimes\cal{A}$ && The adjoint action && $
\ad(a)=a^{(2)}\otimes\k(a^{(1)})a^{(3)}$&\cr\ika
&$\vv\adj \colon \Gamma_{\inv}^\sstar\rightarrow\Gamma_{\inv}^\sstar
\otimes \cal{A}$&& The adjoint action &&
$\adj\pi=(\pi\otimes  \id) \ad$&\cr}
\hrule}
\end{document}